\documentclass[11pt]{amsart}

\usepackage{graphics}
\usepackage{amsthm}
\usepackage{amssymb}
\usepackage{amsmath}
\usepackage{epic}
\usepackage{curves}
\usepackage{epsfig}
\usepackage{times}
\usepackage{amsfonts} 
\input{psfig.sty}

\theoremstyle{plain}
\newtheorem{theorem}{Theorem}[section]
\newtheorem{lemma}[theorem]{Lemma}

\newtheorem{proposition}[theorem]{Proposition}

\theoremstyle{definition}
\newtheorem{definition}[theorem]{Definition}

\newtheorem{example}[theorem]{Example}

\def\action{\bullet}

\def\del{\partial}
\def\ideal#1{\langle #1\rangle}

\newcommand{\baseRing}[1]{\ensuremath{\mathbb{#1}}}
\newcommand{\Z}{\baseRing{Z}}
\def\ZZ{{\mathbb Z}}
\newcommand{\C}{\baseRing{C}}
\def\CC{{\mathbb C}}
\newcommand{\N}{\baseRing{N}}
\def\NN{{\mathbb N}}

\def\RR{{\mathbb R}}

\newcommand{\nin}{\noindent}
\newcommand{\ep}{\medskip}

\def\conv{\operatorname{conv}}

\def\rk{\operatorname{rank}}

\def\vol{\operatorname{vol}}
\def\deg{\operatorname{deg}}

\title[Arbitrary rank jumps for $A$-hypergeometric systems through Laurent 
polynomials]{Arbitrary rank jumps for $A$-hypergeometric systems\\ through
  Laurent polynomials}

\author{Laura Felicia Matusevich}
\address[Laura Felicia Matusevich]{Department of Mathematics \\
  Harvard University \\  
Cambridge, MA}
\email{laura@math.harvard.edu}

\author{Uli Walther}
\address[Uli Walther]{Department of Mathematics, Purdue University,
  West Lafayette,  
IN} 
\email{walther@math.purdue.edu}
\thanks{The first author was partially supported by an NSF
  Postdoctoral Fellowship.
The second author was partially supported by the NSF, the DfG
  and the Humboldt foundation}

\date{March 25, 2004}

\begin{document}

\begin{abstract}
We investigate the solution space of hypergeometric systems of
differential equations in the sense of Gelfand, Graev, Kapranov and
Zelevinsky.  
For any integer $d \geq 2$ we construct a matrix $A_d \in \N^{d \times 2d}$ 
and a parameter vector $\beta_d$ such that the holonomic rank of 
the $A$-hypergeometric system $H_{A_d}(\beta_d)$ exceeds the simplicial
volume $\vol(A_d)$ by at
least $d-1$. The largest previously known gap between rank and volume was 
two.

\ep

\nin Our argument is elementary in that it uses only linear algebra, and
our construction 
gives evidence to the general observation that
rank-jumps seem to go hand in hand with the existence of multiple
Laurent (or
Puiseux) 
polynomial solutions.
\end{abstract}

\maketitle

\section{Introduction}
\label{sec:intro}

\nin A power series $\sum_{t=1}^\infty a(t) x^t$ is geometric, if the
assignment $t\mapsto a(t+1)/a(t)$ is a constant function on $\NN$. 
If the value of these quotients  is always $\lambda$, then clearly
$a(t)=c\cdot \lambda^t$ for some constant $c$.
A natural generalization are the {\em hypergeometric series} for which
$a(t+1)/a(t)$ is a rational function in $t$. The study of such
objects goes back at least to Euler. Gau{\ss} continued this work and 
Kummer and Riemann pioneered the idea of investigating
the differential equations that are satisfied by a given
hypergeometric series.

\ep

\nin Hypergeometric differential equations and their solutions,
hypergeometric functions, are  
a fascinating mixture of algebra, analysis and combinatorics, and 
among the most ubiquitous mathematical objects. They seem to occur naturally 
almost everywhere --- following are just a few examples to illustrate this.
If you try to solve the Laplace partial differential equation by separation
of variables, the Bessel equation appears naturally: its solutions are
hypergeometric \cite{Sobolev}.  
When parameterizing elliptic curves, one encounters theta functions,
which are hypergeometric \cite{Yoshida}. Perhaps
one is  trying to solve a 
polynomial equation 
of degree $n$ in terms of the coefficients: radicals will
not be enough to do  
this if $n>4$, but hypergeometric functions will \cite{Sturmfels-polysol}. 
Or maybe you want to do
least squares approximations on sets of data, and the polynomial basis
you need to  
use involves orthogonal polynomials; all interesting such bases
consist of 
hypergeometric elements \cite{Koekoek}. 
In mirror symmetry, the periods of certain
natural differentials in families of Calabi--Yau toric hypersurfaces satisfy
hypergeometric equations \cite{Cox-Katz}.
If you want to count combinatorial objects and your quantities
satisfy recursions, then this often forces their generating function
to be hypergeometric.  In a recent instance of this phenomenon
involving algebraic geometry, the generating functions for 
intersection numbers on moduli spaces of curves turn
out to be $A$-hypergeometric in the sense of Gelfand, Graev, Kapranov and
Zelevinsky \cite{okounkov}. It is  
this $A$-hypergeometric approach that we shall follow in this article.

\ep

\nin Gelfand, Graev and Zelevinsky defined $A$-hypergeometric systems in the 
mid-eighties, and they were further developed by Gelfand, Kapranov and
Zelevinsky 
(see \cite{GGZ87,GKZ89,GKZ-correction}). 
Before we give the general definition of $A$-hypergeometric systems,
let us consider one example.
\ep

\begin{example}
\label{ex-210}
Let $A$ be the matrix
$\left(\begin{array}{ccc}1&1&1\\2&1&0\end{array}\right)$. 
We consider the integral kernel $\ker_{\ZZ}(A)$ of $A$ consisting of all
$u\in\ZZ^3$ with $A\cdot u=0$.
For our $A$ we have that $\ker_{\ZZ}(A)$ is generated by $u=(1,-2,1)$. We use this vector 
to form the
  operator $\Delta(u)=\frac{\del}{\del x_1}\,\frac{\del}{\del x_3}-
\frac{\del^2}{{\del x_2}^2}$
  by separating the positive part $u_+=(1,0,1)$ from the negative part
  $u_-=(0,2,0)$ of $u$ and then using
  the entries as exponents over the corresponding derivations. 

\ep

\nin From the two rows of the matrix we create the operators 
\begin{eqnarray*}
E_1&=&1\cdot x_1\frac{\del}{\del x_1}+1\cdot x_2\frac{\del}{\del
  x_2}+1\cdot x_3\frac{\del}{\del x_3},\\
E_2&=&2\cdot x_1\frac{\del}{\del x_1}+1\cdot x_2\frac{\del}{\del
  x_2}+0\cdot x_3\frac{\del}{\del x_3}. 
\end{eqnarray*}

\ep

\nin For any pair $\beta=(\beta_1,\beta_2)$ of complex numbers, the
$A$-hypergeometric system is the system of linear partial differential
equations
\begin{eqnarray}
\label{eqn:210}
E_1\action(\varphi)&=&\beta_1\cdot \varphi,\nonumber\\
E_2\action(\varphi)&=&\beta_2\cdot \varphi,\\
\left(\frac{\del^2}{\del x_1\del x_3}-\frac{\del^2}{{\del
  x_2}^2}\right)\action(\varphi)&=&0\nonumber
\end{eqnarray}
where $\varphi$ is a function in the three variables $x_1,x_2,x_3$. One
may interpret $(\beta_1,\beta_2)$ as a multi-degree of the solution
$\varphi$ as we explain now.
First notice that:
\[ 
\left(x_i\frac{\del}{\del x_i}\right) \bullet (x_1{}^{\alpha_1}x_2{}^{\alpha_2}x_3{}^{\alpha_3}) 
= \alpha_i  x_1{}^{\alpha_1}x_2{}^{\alpha_2} x_3{}^{\alpha_3}\;, \quad
i=1,2,3.
\]
This means, using linearity, that for a power series $\varphi(x_1,x_2,x_3) = 
\sum_{\alpha} c_{\alpha} x_1{}^{\alpha_1}x_2{}^{\alpha_2}x_3{}^{\alpha_3}$ we have:
\begin{eqnarray*}
(E_1-\beta_1) \bullet \varphi &=&
\sum_{\alpha} c_{\alpha} (E_1 - \beta_1) \bullet  
\big( x_1{}^{\alpha_1}x_2{}^{\alpha_2}x_3{}^{\alpha_3} \big) \\
&=&\sum_{\alpha} c_{\alpha} \left(x_1 \frac{\del}{\del x_1}+  x_2
 \frac{\del}{\del x_2}+ x_3 \frac{\del}{\del x_3}-\beta_1\right) \bullet  
\big( x_1{}^{\alpha_1}x_2{}^{\alpha_2}x_3{}^{\alpha_3} \big) \\
&=& \sum_{\alpha} c_\alpha\,(\alpha_1 + \alpha_2 + \alpha_3 - \beta_1) \,
 x_1{}^{\alpha_1}x_2{}^{\alpha_2}x_3{}^{\alpha_3}.
\end{eqnarray*}
Thus, if $(E_1 - \beta_1) \bullet \varphi = 0$, then the exponents $\alpha$ appearing
in $\varphi=\sum_{\alpha} c_{\alpha} x_1{}^{\alpha_1}x_2{}^{\alpha_2}x_3{}^{\alpha_3}$
must satisfy:
\[ 
\left[c_{\alpha} \neq 0\right] \Longrightarrow \left[\alpha_1 + \alpha_2 + \alpha_3 =
\beta_1\right]. 
\]
A similar computation using $E_2$ instead of $E_1$ yields:
\[
\left[c_{\alpha} \neq 0\right] \Longrightarrow \left[2\alpha_1 + \alpha_2  = \beta_2\right]
\]
and the two implications combine to
\begin{equation}
\label{eqn:homog}
\left[c_{\alpha} \neq 0\right] \Longrightarrow \left[A \cdot \alpha = \beta\right]. 
\end{equation}
Let us define the {\em multi-degree} of $x_i$ to be the $i$th column of $A$:
\[ \deg(x_1) = \left( \begin{array}{c} 1 \\ 2 \end{array} \right); \;
 \deg(x_2) = \left( \begin{array}{c} 1 \\ 1 \end{array} \right); \;
 \deg(x_3) = \left( \begin{array}{c} 1 \\ 0 \end{array} \right); \]
hence the multi-degree of a monomial is given by:
\[ \deg(x_1{}^{\alpha_1}x_2{}^{\alpha_2}x_3{}^{\alpha_3}) = A \cdot \alpha.\]
Now equation (\ref{eqn:homog}) translates into:
\begin{center}
If $\varphi=\sum_{\alpha} c_{\alpha} 
x_1{}^{\alpha_1}x_2{}^{\alpha_2}x_3{}^{\alpha_3}$ 
is killed by $E_1-\beta_1$ and $E_2 -\beta_2$,
then

\ep

$\left[ c_{\alpha} \neq 0 \right] \Longrightarrow
\left[\deg(x_1{}^{\alpha_1}x_2{}^{\alpha_2}x_3{}^{\alpha_3}) = \beta \right]$.
\end{center}

\bigskip

\nin To illustrate one point made in the introduction above, let
$\beta=(0,-1)$. It is well-known and easy to verify that then the two roots
$z_{1,2}=\frac{-x_2\pm\sqrt{{x_2}^2-4x_1x_3}}{2x_1}$ 
of the polynomial $x_1z^2+x_2z+x_3$ in the variable $z$ with
indeterminate coefficients $x_1,x_2,x_3$ are solutions of the system
(\ref{eqn:210}). In turn, one can use the system of partial
differential equations to obtain a formula of the roots as a
hypergeometric series:
\[
z_{1,2}=
\frac{-x_2}{2x_1}\pm\left(\frac{x_2}{2x_1}-
\frac{x_3}{x_2}\sum_{t=0}^\infty
\frac{1}{t+1}{2t\choose
  t}\left(\frac{x_1x_3}{{x_2}^{2}}\right)^t\right).
\]
\end{example}

\ep\ep

\nin We now come to the definition of a general $A$-hypergeometric
system. 
We begin with taking 
an integer $d\times n$ matrix $A=(a_{i,j})$ of full rank $d$ and
a complex parameter vector $\beta$. As in the example
we form for $1\le i\le d$ the operators 
\[
E_i=\sum_{j=1}^n a_{i,j}\, x_j \frac{\partial}{\partial x_j}
\]
from the rows of $A$.

\begin{definition}
\label{def:A-hyp}
The {\em $A$-hypergeometric system with parameter $\beta$}, denoted
$H_A(\beta)$, 
is the following system of linear partial 
differential equations with polynomial coefficients for the function
$\varphi=\varphi(x_1,\ldots,x_n)$:
\begin{align*}
{\rule[-3ex]{0cm}{6ex}}
E_i \action(\varphi)  &  = \beta_i \cdot \varphi &\quad \quad i=1,\dots d;\\ 
\phantom{xxxxxxxxxxxxxxxxxx}\left( \prod_{u_i>0}
\frac{\partial^{u_i}}{\partial {x_i}^{u_i}} \right) \action(\varphi) & = 
\left( \prod_{u_i<0} \frac{\partial^{-u_i}}{\partial {x_i}^{-u_i}} \right) 
\action(\varphi)  &
\quad   \quad 
\mbox{for all} \; u  \in \ker_{\ZZ}(A).
\end{align*}
The first $d$ equations above are called {\em homogeneity
  conditions}, the
remaining equations are called {\em toric equations}.
\end{definition}

\ep 

\nin For notational convenience we shall from now on abbreviate the
derivation $\frac{\del}{\del x_i}$ by simply $\del_i$. 
Then
$R_A=\CC[\del_1,\ldots,\del_n]$ is the ring of $\CC$-linear 
differential operators
with constant coefficients.
Let us view Example \ref{ex-210} in the light of 
our definition of general hypergeometric systems. 
In Definition \ref{def:A-hyp} there are infinitely many
toric equations, one for each element $u$ of $\ker_\ZZ(A)$. On the
other hand, 
in (\ref{eqn:210}) we listed only one such,
$\Delta(u)\action\varphi=0$ with
$u=(1,-2,1)$. 
Yet it turns out that no
information is lost. Namely, if $A$ is the matrix of Example
\ref{ex-210} and $v\in\ker_\ZZ(A)$ then up to sign $v=(k,-2k,k)$ for
some natural number $k$. It follows that, again up to sign, 
\begin{eqnarray*}
\Delta(v)&=&({\del_1\del_3})^k-{\del_2}^{2k}\\
&=&\left((\del_1\del_3)^{k-1}+(\del_1\del_3)^{k-2}{\del_2}^{2}+
 (\del_1\del_3)^{k-3}{\del_2}^{4}+\cdots+{\del_2}^{2k-2}\right)\cdot
 \left(\del_1\del_3-{\del_2}^2\right).
\end{eqnarray*}
So
if $\varphi$ is annihilated by $\Delta(u)$ then it is also annihilated by
$\Delta(v)$ for all other $v\in\ker_\ZZ(A)$.

\ep

\nin More generally, 
it turns out that for any matrix $A$ one always only needs
to look at a 
finite number of toric equations; in order to explain the reasons for
this we simplify our notation a bit as follows. 
In the remainder of the paper we would like to use 
multi-index notation: if $u\in\ZZ^n$ we mean by $x^u$ 
the (Laurent) monomial ${x_1}^{u_1}{x_2}^{u_2}\cdots {x_n}^{u_n}$; a
similar convention shall be used for $\del^u$. 
Also, if $u \in \ZZ^n$, we write $u=u_+ - u_-$, where:
\[ 
(u_+)_i = \max\{u_i,0\}\, , \quad (u_-)_i = \max\{-u_i,0\}.
\] 
With this notation, the toric operator $\Delta(u)=\prod_{u_i>0}
\frac{\partial^{u_i}}{\partial {x_i}^{u_i}} - \prod_{u_i<0}
  \frac{\partial^{-u_i}}{\partial 
  {x_i}^{-u_i}}$ in $H_A(\beta)$ corresponding 
to $u\in\ker_\ZZ(A)$
becomes $\del^{u_+}-\del^{u_-}$. 
Let $I_A$ be the {\em toric ideal} in
$R_A$ generated 
by all $\Delta(u)=\del^{u_+}-\del^{u_-}$ with $u\in\ker_\ZZ(A)$. Since
$R_A$ is Noetherian, there is a finite set of generators for this
ideal. In fact, since $I_A$ is generated by {\em binomials}, this finite
generating set will consist of binomials and hence be of the form
$\{\Delta(v_1),\ldots,\Delta(v_k)\}$ for some elements
$v_1,\ldots,v_k$ in $\ker_\ZZ(A)$. Indeed, there are simple algorithms
to find such a collection $\{v_i\}_{i=1}^k$, see \cite{Sturmfels}.

\ep

\nin Although we will not use this, we would like to mention that
by a theorem of
Stafford \cite{Stafford} the entire $A$-hypergeometric system is
equivalent to a linear system of just {\em two} differential
equations. However, these two equations are very complicated 
since they have to carry a lot of information.

\ep

\nin Since $H_A(\beta)$ is a {\em linear} system of equations, the
set of its holomorphic solutions on a simply connected open set in
$\CC^n$  forms a vector space over the
complex numbers. The dimension of this vector space we 
shall call the {\em rank} of $H_A(\beta)$ and denote it
by $\rk(H_A(\beta))$. 
Somewhat surprisingly, the rank turns out to be finite for any
choice of $A$ and $\beta$ --- this is a highly unusual event for systems
partial differential equations.

\ep

\nin So one of the most basic questions one might ask about the
$A$-hypergeometric system 
$H_A(\beta)$ is:

\bigskip

\nin {\bf Question A:} What is the rank of $H_A(\beta)$?

\bigskip

\nin A first answer to this question was given by Gelfand, Kapranov and Zelevinsky 
\cite{GKZ89,GKZ-correction} who found that under a certain condition
on the ideal $I_A$ called {\em Cohen--Macaulayness}, 
$\rk(H_A(\beta))$ is actually {\em independent} of $\beta$. To
describe this condition, consider the polynomial ring
$R_A=\CC[\del_1,\dots,\del_n]$ from above and its quotient $S_A=R_A/I_A$. 
Then one calls $I_A$ {\em Cohen--Macaulay} if and only if there are
$d=\rk(A)$ linear forms $L_1,\ldots,L_d$ in $R_A$ such that for all
$1\le i\le d$ the form $L_i$ is a non-zerodivisor on
$S_A/\ideal{L_1,\ldots,L_{i-1}}$. This property is a way of
allowing singularities to occur in $S_A$ while preserving many good
algebraic properties. By a theorem of Hochster \cite{Hochster}, one
particular class of 
Cohen--Macaulay examples is provided by those matrices $A$ for which
the collection $\NN A$ of all $\NN$-linear combinations of the columns
of $A$ is {\em saturated}. This condition means that if a lattice point 
$p\in\ZZ^d$ has some
multiple $p+\cdots +p$ in $\NN A$, then 
$p$ itself is already in $\NN A$. Such saturated semigroups arise
naturally as
the 
collection of all lattice points inside the {\em positive cone} $\RR_+
v_1+\cdots +\RR_+ v_k$ of $k$ lattice points $v_1,\ldots,v_k \in
\ZZ^{d}$. Our Example \ref{ex-210} is of this type with 
$v_1=\left(\begin{array}{c}1\\2\end{array}\right)$ and 
$v_2=\left(\begin{array}{c}1\\0\end{array}\right)$.

\ep

\nin Under the assumption of Cohen--Macaulayness,  
a completely explicit  combinatorial formula for the rank was provided in
\cite{GKZ89,Adolphson-Duke94}. Let us describe this formula.
Form a polytope $Q_0$ 
by taking the convex hull of the columns of $A$ and the origin, 
pictured as points in $\RR^d$.
Since $A$ has full rank
this polytope
has dimension $d$. Then the {\em simplicial} or {\em normalized volume} of $A$, denoted 
by $\vol(A)$ equals the product of $d!$ and the usual Euclidean volume
of $Q_0$ (so that, for example, a standard $d$-simplex has simplicial
volume equal to $1$).
With this notation, if $I_A$ is Cohen--Macaulay, then 
the rank $\rk(H_A(\beta))$ of the hypergeometric
system to $A$ and $\beta$ agrees with the simplicial volume $\vol(A)$
no matter what the parameter $\beta \in \C^d$ is.

\ep

\nin Several authors have expanded on these results, usually
in the {\em homogeneous} case where 
all the columns of $A$, considered as points
in $\RR^d$, lie in a hyperplane not containing  the origin. For example,  
Adolphson \cite{Adolphson-Duke94}
showed that even if $A$ fails to be Cohen--Macaulay then the formula
$\rk(H_A(\beta))=\vol(A)$ 
is valid for
{\em almost every} $\beta$. 
If $A$ is homogeneous, but under no other conditions 
on either $A$ or $\beta$, we always have $\rk(H_A(\beta))\geq \vol(A)$ as
was shown by Saito, Sturmfels and Takayama \cite{SST}.
Considering these results, the natural question is: 

\bigskip

\nin {\bf Question B:} Are there actually any examples where $\rk(H_A(\beta)) > \vol(A)$ ?

\bigskip

\nin The answer is ``yes'', and the first and smallest example of this type 
was given in 
\cite{Sturmfels-Takayama}; we will revisit  it in 
Example \ref{example:0134}. Experimental studies showed that
constructing rank-jumping examples $(A,\beta)$ is very hard since they
are quite rare; this accounts for the 10-year delay between the first
results on $A$-hypergeometric functions and the discovery of the first
rank-jump.

\ep

\nin One reason that makes rank-jumps very interesting is that they seem to
coincide with the existence of very nice solutions: contrary to typical
solutions which are proper power series, in all cases that are known
to the authors  the ``extra''
solutions at a rank-jump are 
{\em Laurent  polynomials} (or Puiseux polynomials, if the exponents
are non-integral); this
fact is not well understood yet. 
Viewing the results of \cite{Adolphson-Duke94,GKZ89,SST} in the light of
Example \ref{example:0134},
 one is then lead to three more precise questions:

\bigskip

\nin {\bf Questions C:}
\begin{enumerate}
\item Which matrices $A$ allow for rank jumps?
\item If $A$ has a rank jump at all, which parameters are rank-jumping?
\item If $\beta$ is a rank-jumping parameter for $A$, by how much does
  the rank exceed the volume?
\end{enumerate}

\ep

\nin The first two questions have been recently answered in full
\cite{MMW}.
In the present article we are interested in the third question and
investigate the  
possible magnitude of the gap between rank and volume. 
There is a known upper bound for the rank in terms of the volume given by
$\rk(H_A(\beta)) \leq 2^{2d} \cdot\vol(A)$, see  \cite[Corollary
  4.1.2]{SST}. 
It is believed that this exponential upper bound
is not optimal. In fact, 
until now no
example had been known in which the rank exceeds the volume by three
or more. 

\ep

\nin The goal of this article to describe a family of examples that exhibit
arbitrarily large 
rank jumps, we shall prove:

\begin{theorem}
\label{theorem:main}
For any $d \in \Z_{>1}$ there exists a $d\times(2d)$-matrix 
$A_d$ and a parameter $\beta_d\in\CC^d$ such that
\[ \rk(H_{A_d}(\beta_d))-\vol(A_d) \geq d-1. \]
\end{theorem}

\ep

\nin In contrast to the substantial amount of algebra and analysis
that is needed 
to prove most of the results quoted above, the proof of our 
result is completely elementary, requires only a knowledge
of linear algebra and is based on constructing Laurent polynomial solutions. 


\section{The first rank-jump example}

\nin We now present 
a major player in
our later constructions:
the first ever rank-jumping example.

\ep

\begin{example}
\label{example:0134}
Let $\beta=(\beta_1,\beta_2)$ and 
\[ 
A_2 =\left( \begin{array}{cccc} 1 & 1 & 1 & 1 \\ 0 & 1 & 3 & 4
\end{array} \right).
\]
Then $I_{A_2}$ is generated by 
\begin{eqnarray*} 
\partial_2\partial_3 &-&\partial_1 \partial_4,\\
\partial_1^2\partial_3 &-&\partial_2^3, \\
\partial_2\partial_4^2&-&\partial_3^3,\\ 
\partial_1\partial_3^2&-&\partial_2^2\partial_4
\end{eqnarray*}
and there are two homogeneity conditions:
\begin{eqnarray*}
(x_1 \partial_1 + x_2 \partial_2 + \phantom{3}x_3 \partial_3 +
  \phantom{4}x_4 \partial_4 -\beta_1)\action(\varphi)&=&0,\\ 
(x_2\partial_2 + 3x_3 \partial_3 + 4 x_4
  \partial_4-\beta_2)\action(\varphi)&=&0.
\end{eqnarray*}
In this case, 
\[ 
\rk(H_{A_2}(\beta_1,\beta_2)) = \left\{ \begin{array}{ll}
4 = \vol(A_2) \quad & {\rm if\ }\, (\beta_1,\beta_2) \neq (1,2), \\
5 & {\rm if\ }\, (\beta_1,\beta_2)=(1,2). \end{array} \right. \]
\end{example}

\ep

\nin Example \ref{example:0134}
was completely
analyzed in  \cite{Sturmfels-Takayama}. We refer to that article for a proof
that $(1,2)$ is indeed the unique parameter for which rank exceeds
volume. We now present
an explicit basis for the solution space of $H_{A_2}(1,2)$.

\ep

\begin{theorem}[Proposition 4.1 \cite{Sturmfels-Takayama}]
\label{thm:basis}
Let
\[ 
u^{(1)} = (1/2,0,0,1/2),\,\,\, u^{(2)} = (1/4,1,0,1/4),\,\,\, u^{(3)}
=(1/4,0,1,-1/4) 
\]
and put for $i=1,2,3$
\[
\Omega_i=\left\{(a,b)\in\ZZ^2: u_2^{(i)} +4a \geq 3b,\, u_3^{(i)} +b \geq
0\right\}.
\]
Consider for $i=1,2,3$ the functions
\[ f_i = \sum_{(a,b)\in \Omega_i}
c_{a,b}\, x^{u^{(i)}+a(-3,4,0,-1)+b(2,-3,1,0)}
\]
where
\[ 
c_{a,b} = \frac{1}{\Gamma(u_1^{(i)}-3a+2b+1)
\Gamma(u_2^{(i)}+4a-3b+1)
\Gamma(u^{(i)}_3+b+1)
\Gamma(u^{(i)}_4-a+1)} 
\]
and $\Gamma$ denotes the usual gamma function.
If one sets
\[ p_1 = \frac{{x_2}^2}{x_1}, \qquad p_4 = \frac{{x_3}^2}{x_4} \]
then
the five functions $p_1,p_4,f_1,f_2,f_3$ are a basis for the solution space of
$H_{A_2}(1,2)$.\qed
\end{theorem}


\section{Constructing arbitrary jumps}

\nin We are now ready to provide, for given $d \geq 2$, a $d \times 2d$ 
matrix $A_d$ and a parameter $\beta_d \in \N^d$ such that
\[ \rk(H_{A_d}(\beta_d)) \geq \vol(A_d) + d-1 . \]
As we mentioned before,  previously no example existed where the gap
between rank and volume 
exceeds two.

\ep

\nin If $d=2$, Example \ref{example:0134} will do. So for the remainder of
this article we fix an integer $d \geq 3$, 
and we write  
$A$ and $\beta$ instead of $A_d$ and $\beta_d$ in order to simplify notation.

\ep

\nin Let $e_1,\dots ,e_d$ be the standard basis vectors
in $\C^d$. Define $a_1,\dots , a_{2d} \in \N^d$ as follows:
\begin{align*}
a_1 & = (1,0,\dots,0,0), \\
a_2 & = (1,0,\dots,0,1), \\
a_3 & = (1,0,\dots,0,3), \\
a_4 & = (1,0,\dots,0,4), 
\end{align*}
while if $3 \leq k \leq d-1$, set
\[ 
a_{2k-1} = e_1+e_{k-1}, \qquad a_{2k}=e_1+e_{k-1}+e_{d}.
\]
Thus
\[
A= \left( \begin{array}{ccccccccccc}
1      &  1     &  1  &  1  &  1  &  1  &  1  &  1  & \cdots &  1     &  1     \\ 
0      &  0     &  0  &  0  &  1  &  1  &  0  &  0  &        &  0     &  0     \\
0      &  0     &  0  &  0  &  0  &  0  &  1  &  1  &        &  0     &  0     \\
\vdots & \vdots &     &     &     &     &     &     & \ddots & \vdots & \vdots \\ 
0      &  0     &  0  &  0  &  0  &  0  &  0  &  0  &        &  1     &  1     \\
0      &  1     &  3  &  4  &  0  &  1  &  0  &  1  & \cdots &  0     &  1  
\end{array}
\right)
\]

\ep

\nin Now let
\[ 
\beta = (1,0,\dots,0,2).
\]
We shall prove

\ep
 
\begin{theorem}
\label{thm:secondmain}
For the matrix $A$ and parameter $\beta$ introduced above, we have:
\[ \rk(H_A(\beta))-\vol(A) \geq d-1 .\]
\end{theorem}

\ep

\nin We will prove this theorem in a series of lemmas. First we will
compute the simplicial
volume $\vol(A)$; after this is done, we will exhibit the required
number  of linearly
independent solutions of $H_A(\beta)$.

\ep

\begin{lemma}
The simplicial volume of $A$ is $d+2$.
\end{lemma}

\ep

\begin{proof}
Let $Q=\conv(A)$, the convex hull of the columns of $A$. 
Since the columns of $A$ all lie in the hyperplane $t_1=1$ of $\RR^d$,
the convex hull $Q_0$ of the origin and the columns of $A$ form a
pyramid of height one over $Q$. Hence 
the simplicial volume of $Q_0$
is equal to
the simplicial volume of $Q$; we compute the latter. 

\ep

\nin The polytope $Q$ is the union of two others: the prism $P$ (over the
standard $(d-2)$-simplex with vertices  $p_1,p_5,p_7,\ldots,p_{2d-1}$) 
whose vertices are the columns of:
\[
\left( \begin{array}{ccccccccccc}
1 & 1 & 1  &  1  &  1  &  1  & \cdots &  1     &  1     \\ 
0 & 0 & 1  &  1  &  0  &  0  &        &  0     &  0     \\
0 & 0 & 0  &  0  &  1  &  1  &        &  0     &  0     \\
\vdots & &    &     &     &     & \ddots & \vdots & \vdots \\ 
0 & 0 & 0  &  0  &  0  &  0  &        &  1     &  1     \\
0 & 1 & 0  &  1  &  0  &  1  & \cdots &  0     &  1  
\end{array}
\right), 
\]
and the $(d-1)$-simplex $S$ whose vertices $p_2,p_4,p_6,\ldots,p_{2d}$
are the columns of:
\[
\left( \begin{array}{ccccccccccc}
  1     &  1  &  1  & \cdots &  1        \\ 
  0     &  0  &  1  &        &  0        \\
 \vdots &     &     & \ddots & \vdots    \\ 
  0     &  0  &  0  &        &  1        \\
  1     &  4  &  1  & \cdots &  1    
\end{array}
\right).
\]
In Figure \ref{fig:decomp} we see the decomposition of $Q$ into 
the prism $P$ and the simplex $S$ for $d=4$.
\begin{figure}[hbtp]
  \centerline{\psfig{
       figure=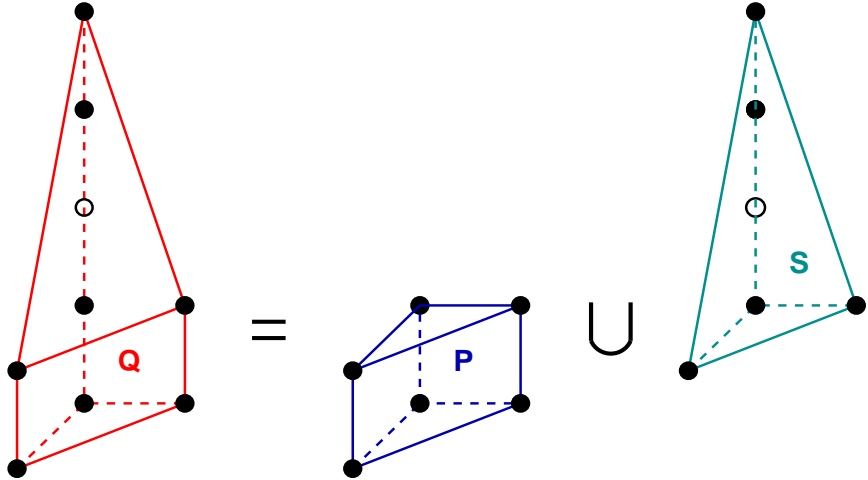,height=2.5in,width=4.5in}}      
    \caption{Decomposing $Q=\conv(A)$ for $d=4$}{\label{fig:decomp}
    }
\end{figure}
Since the prism has height one, its Euclidean volume equals the
Euclidean volume of  
its base, the standard $(d-2)$-simplex with Euclidean volume 
$\frac{1}{(d-2)!}$. Thus, $P$ has simplicial volume
$\frac{(\dim(P))!}{(d-2)!}=d-1$. 

\ep

\nin On the other hand, $S$ is a pyramid of height three over a standard
simplex, 
and so its simplicial volume is $3$. 
This implies that $\vol(Q) = (d-1)+3 = d+2$.
\end{proof}

\ep

\nin The next step in our proof is to construct $2d+1$ 
solutions of $H_A(\beta)$. 
In order to do this we need to understand the integer kernel
of $A$, because the toric equations are constructed directly
from these elements. In particular, we will identify positive and
negative coordinates of certain elements in $\ker_{\Z}(A)$. The other
important ingredient 
is finding integer solutions of $A\cdot u = \beta$. The fact that the
coordinates 
of $\beta$ are small positive integers will facilitate this search.
However, we start with showing 
that any solution of $H_{A_2}(1,2)$ is a solution of our system.

\ep

\begin{lemma}
\label{lemma:prevsols}
Let $\psi$ be a solution of $H_{A_2}(1,2)$. Then $\psi$ is a solution
of $H_A(\beta)$. 
In particular, the functions $p_1$, $p_4$, and $f_1$, $f_2$, $f_3$ 
from Theorem \ref{thm:basis} are linearly independent solutions
of $H_A(\beta)$.
\end{lemma}

\ep

\begin{proof}
It is easy to see that $\psi$ is a solution of the homogeneity equations
\begin{align*}
\phantom{xxxxxxxxxxxxxxxxxxxxxxxx}&\sum_{j=1}^{2d} a_{i,j}\,x_j \partial_j\action(\psi) =
\beta_i\cdot\psi, &\qquad\qquad i=1,\dots ,d. 
\end{align*}
Hence we only need to verify that $\psi$ is annihilated by the toric operators
$\Delta(u)=\partial^{u_+} - \partial^{u_-}$ for all
$u_+-u_- =u\in \ker_{\ZZ}(A)$.
We now
study the integer kernel $A$.
Since $A$ is of full rank $d$ and the
columns of the following $(2d\times d)$-matrix $B$ are linearly
independent, the columns of $B$ form a basis for the
kernel of $B$ over the rational numbers:
\[ 
B = \left( 
\begin{array}{rrrrrrr}
 1 &  1 &  1 &  1 &  1 & \cdots &  1 \\
-2 & -1 & -1 & -1 & -1 & \cdots & -1 \\
 2 & -1 &  0 &  0 &  0 & \cdots &  0 \\
-1 &  1 &  0 &  0 &  0 & \cdots &  0 \\
 0 &  0 & -1 &  0 &  0 & \cdots &  0 \\
 0 &  0 &  1 &  0 &  0 & \cdots &  0 \\
 0 &  0 &  0 & -1 &  0 & \cdots &  0 \\
 0 &  0 &  0 &  1 &  0 & \cdots &  0 \\
 0 &  0 &  0 &  0 & -1 & \cdots &  0 \\
 0 &  0 &  0 &  0 &  1 & \cdots &  0 \\
\vdots & & & & \vdots &   & \vdots \\
 0 &  0 &  0 &  0 &  0 & \cdots & -1 \\
 0 &  0 &  0 &  0 &  0 & \cdots &  1 \\
\end{array}
\right)
\]
Using rows $1,2,5,7,9,\ldots,2d-1$ we see that
the greatest common divisor of the maximal minors of $B$ is $1$. This
implies that the columns of
$B$ are actually 
a basis for the integer kernel $\ker_\ZZ(A)$: any element
of $\ker_{\Z}(A)$ is an {\em integer} linear combination of
the columns of $B$.

\ep

\nin Choose a toric operator $\partial^{u_+}-\partial^{u_-}$ where $u =
B\cdot z$ for 
some $z \in \Z^d$. If $z_i \neq 0$ for some $i\geq 3$, then $u_{2i-1}$ and
$u_{2i}$ will be nonzero, with opposite signs. This means that one of
the monomials 
in $\partial^{u_+}-\partial^{u_-}$ will contain $\partial_{2i-1}$, and
the other 
will contain $\partial_{2i}$. Since $\psi$ does not contain the variables
$x_{2i-1}$ nor $x_{2i}$, it follows that both monomials annihilate $\psi$ and 
therefore $( \partial^{u_+}-\partial^{u_-})\action( \psi) = 0$.

\ep

\nin It remains to consider the case when only $z_1$ and $z_2$ are allowed
to be nonzero. 
But in that case $u=B\cdot z$ gives a 
toric operator inside $H_{A_2}(1,2)$, and $\psi$ was assumed to be a
solution of that system.
\end{proof}

\ep

\nin We note that there are no polynomial solutions for $H_A(\beta)$ since
any such solution would have to have multi-degree $\beta$ and $\beta$
is not an $\NN$-linear combination of the multi-degrees of the $x_i$, which
are the columns of $A$.
We will now construct Laurent polynomial solutions for $H_A(\beta)$,
one for each vertex 
of the polyhedron $Q=\conv(A)$. These vertices 
are the columns $a_1$, $a_4$, $a_5,\dots ,a_{2d}$
of the matrix $A$. The correspondence between the Laurent polynomials $p_i$
and the vertices $a_i$ will be given by 

\ep

\begin{center}
$p$ is associated to $a_i$ if no variable but $x_i$ occurs
  in any denominator of $p$.
\end{center}

\ep

\nin For the
vertices $a_1$ and $a_4$ we already have such solutions,
namely the Laurent monomials
$p_1={x_2}^2/x_1$ and $p_4={x_3}^2/x_4$. 
So we need to construct Laurent polynomial solutions $p_i$ of
$H_A(\beta)$ 
associated to $a_5,\ldots,a_{2d}$, and since $H_A(\beta)$ does not
have polynomial solutions, these are proper fractions.

\ep

\nin If $p = \sum c_{\alpha} x^{\alpha}$ is a Laurent polynomial solution
of $H_A(\beta)$ then 
the homogeneity equations imply that $A\cdot \alpha = \beta$ for any
$\alpha$ such that 
$c_{\alpha} \neq 0$. Hence the possible Laurent monomials appearing in
a Laurent solution $p_i$
of $H_A(\beta)$ associated to $a_i$ are of the
form $x^{\alpha}$ where $A\cdot \alpha = \beta$, $\alpha \in \Z^n$ and 
only $\alpha_i$ is a negative integer.

\ep

\nin Let us search for all such vectors $\alpha$ when $i=5$. 
Since the second coordinate of $\beta$
is zero and only the columns $a_5$ and $a_6$
of $A$ have nonzero second coordinates, 
we must have $\alpha_6 = -\alpha_5>0$.

\ep

\nin Then
\[ 
\alpha_5 a_5 + \alpha_6 a_6 = \alpha_6 e_d. 
\]
Note that $A$ has no negative entries.
As $\alpha_i\geq 0$ for $i\not =5$,
$A\cdot\alpha=(1,0,\ldots,0,2)$  is in each component
bounded from below  by $\alpha_5 a_5 + \alpha_6 a_6$,
so $\alpha_6$ equals $1$ or $2$.
Moreover, every $a_i$ has a $1$ in the first coordinate and so
$\alpha$ has precisely one more
nonzero entry besides $\alpha_5$ and $\alpha_6$; this entry
will be a $1$.
Now if $\alpha_j =1$ for any $j>6$  then $A\cdot \alpha$ will have a $1$
in a place where $\beta$ has a zero. Therefore, the third nonzero coordinate
of $\alpha$ must be one of $\alpha_1$, $\alpha_2$, $\alpha_3$ or $\alpha_4$.
If $\alpha_6=1$, we get $\alpha = (0,1,0,0,-1,1,0,0,\dots , 0)$ while
for $\alpha_6=2$ we get $\alpha = (1,0,0,0,-2,2,0,0,\dots,0)$; there
is no other choice.

\ep

\nin This gives us two possible monomials to make a Laurent polynomial solution of 
$H_A(\beta)$ where only $x_5$ is in the denominator, namely the monomials
$\frac{x_2x_6}{x_5}$ and $\frac{x_1{x_6}^2}{{x_5}^2}$.
Neither of these Laurent monomials is a solution for $H_A(\beta)$, but
a suitable linear combination is: 

\ep

\begin{lemma}
\label{lemma:p5}
The function 
\[p_5 = \frac{x_2x_6}{x_5} - \frac{1}{2} \frac{x_1{x_6}^2}{{x_5}^2} \]
is a solution of $H_A(\beta)$.
\end{lemma}

\ep

\begin{proof}
By our construction,  $p_5$ is a solution of the homogeneity equations, 
\begin{align*}
\phantom{xxxxxxxxxxxxxxxxxxxxxxxx}&\sum_{j=1}^{2d} a_{i,j}\,x_j \partial_j\action(\psi) =
\beta_i\cdot\psi, &\qquad\qquad i=1,\dots ,d, 
\end{align*}
because the 
exponents appearing in it satisfy $A\cdot \alpha = \beta$. Now we need
to see that $p_5$ is a solution to 
\[ 
\left(\partial^{u_+} - \partial^{u_-}\right)\action(p_5)=0 
\]
whenever
$u_+-u_- =u\in \ker_{\Z}(A)$.

\ep

\nin Recall that $\ker_{\Z}(A)$ has a $\Z$-basis consisting
of the columns of the matrix $B$.
Let us look at a toric equation $\partial^{u_+} - \partial^{u_-}$, where
$u = B\cdot z$ for some integer vector $z\in \Z^d$.
If $z_i \neq 0$ for some $i>3$, then $u_{2i-1}$ and $u_{2i}$ are nonzero with opposite
signs. Then $\partial_{2i-1}$ and $\partial_{2i}$ appear in different monomials in
$\partial^{u_+}-\partial^{u_-}$ while $p_5$ does not contain either 
of the variables $x_{2i-1}$ or $x_{2i}$. This means that 
\[ 
(\partial^{u_+} - \partial^{u_-})\action( p_5) = 0 \quad \mbox{for} \; u_+-u_- = B\cdot z
\; \mbox{with} \; z_i \neq 0 \; \mbox{for some} \; i>3.
\]
So let us now look at $u =B \cdot z$ for $z$ such that $z_i=0$, $i=4,5,\dots, d$.
Then the only (possibly) nonzero coordinates of $u$ are the following:
\begin{align*}
u_1 &=  z_1+z_2+z_3, \\
u_2 &=  -2z_1-z_2-z_3, \\
u_3 &=  2z_1-z_2, \\
u_4 &=  -z_1+z_2, \\
u_5 &=  -z_3, \\
u_6 &=  z_3, \\
\end{align*}
with all $z_i\in \ZZ$.
If $u_3$ and $u_4$ are both nonzero and have different signs, then the
fact that 
$p_5$ contains neither $x_3$ nor $x_4$ implies that 
$(\partial^{u_+} - \partial^{u_-})\action( p_5) = 0$. This means that we need to study
three cases:
\begin{enumerate}
\item $u_3=u_4=0$,
\item $0\le u_3, u_4$ and not both $u_3$ and $u_4$ vanish,
\item $0\geq u_3, u_4$ and not both $u_3$ and $u_4$ vanish.
\end{enumerate}
In Case (1), we have $z_1=z_2=0$. If $|z_3|\geq 2$, then we have
$\partial_1^2$ and 
$\partial_2^2$ in different monomials of $\partial^{u_+}-\partial^{u_-}$, which
implies that $(\partial^{u_+}-\partial^{u_-})\action(p_5)=0$. In the
remaining case $|z_3|=1$ one finds
\[ 
(\partial_1\partial_6-\partial_2\partial_5)\action( p_5) 
= 0-\frac{1}{2}\frac{2x_6}{{x_5}^2} -
\frac{-x_6}{{x_5}^2} -0= 0.
\]

\ep 

\nin In Case (2) one sees immediately that $\del^{u_+}$ kills $p_5$ since
$u_3$ or $u_4$ will be positive and $p_5$ does not involve either
variable. So we need to show that $\del^{u_-}$ also kills $p_5$.
From the given inequalities one
deduces that either $z_1=z_2=1$ or that $z_1\geq 1$
and $z_2\geq 2$. In the latter situation $u_2\le
-4-z_3\le -2$, so $\del^{u_-}$ contains ${\del_2}^2$ and hence kills
$p_5$. We now consider the case  $z_1=z_2=1$.
Clearly if $z_3<-2$ then $\del^{u_-}$ contains ${\del_6}^3$ and hence
kills $p_5$. 
If $z_3=-2$ then $\del^{u_-}=\del_2{\del_6}^2$ kills
$p_5$. Finally, if $z_3\geq -1$ then $\del^{u_-}$ contains ${\del_2}^2$
and kills $p_5$.

\ep

\nin Case (3) is entirely parallel to Case (2), with signs reversed.
\end{proof}

\ep

\nin The construction of $p_6$ goes along the same lines as the construction of 
$p_5$. First
we find that the only solutions of $A\cdot \alpha =\beta$ with $\alpha_i \in \Z_{\geq 0}$
for $i \neq 6$ and $\alpha_6 \in \Z_{<0}$ are the vectors
\[ 
(0,0,1,0,1,-1,0,\dots,0) \quad \mbox{and} \quad
(0,0,0,1,2,-2,0,\dots,0) ,
\]
Then we propose 
\[
p_6 = \frac{x_3x_5}{x_6}-\frac{1}{2}\frac{x_4{x_5}^2}{{x_6}^2}. 
\]
A similar analysis as in Lemma \ref{lemma:p5} shows that, except for 
$\partial_3\partial_6-\partial_4\partial_5$, every generator
$\partial^{u_+}-\partial^{u_-}$ of $I_A$ has the property that
both $\partial^{u_+}$ and $\partial^{u_-}$ annihilate $p_6$. Now to
establish $p_6$ as solution of $H_A(\beta)$ reduces to
 checking that
\[ (\partial_3\partial_6-\partial_4\partial_5)\action(p_6) = 
\frac{-x_5}{{x_6}^2} - \frac{-1}{2}
\frac{2x_5}{{x_6}^2} = 0.\]

\ep

\nin More generally, adapting the notation, we obtain:

\ep

\begin{proposition}
\label{propo:pi}
The two functions
\[ p_{2i-1} = \frac{x_2x_{2i}}{x_{2i-1}} - 
\frac{1}{2}\frac{x_1{x_{2i}}^2}{{x_{2i-1}}^2}, \qquad   
   p_{2i}   = \frac{x_3x_{2i-1}}{x_{2i}} -
   \frac{1}{2}\frac{x_4{x_{2i-1}}^2}{{x_{2i}}^2}
\] 
are solutions of $H_A(\beta)$ for every integer $i$ with $3 \leq i\leq d$.
\end{proposition}

\ep

\nin We can now complete the proof of our main result.

\ep

\begin{proof}[Proof of Theorem \ref{thm:secondmain}]
The functions $f_1,f_2,f_3$ and $p_1,p_4,p_5,\dots ,p_{2d}$ are $2d+1$ 
solutions of $H_A(\beta)$. By Theorem \ref{thm:basis}, the first five are
linearly independent. Since for $i>4$ the Laurent solution $p_i$ has a
pole in $x_i$ and since $x_i$ does not occur in the solutions
$f_1,f_1,f_3,p_1,p_4,\ldots,p_{i-1}$ we conclude that all these
solutions are  linearly independent. It
follows that  $\rk(H_A(\beta)) \geq 2d+1$.
Using $\vol(A) = d+2$, we conclude that
\[ 
\rk(H_A(\beta)) - \vol(A) \geq d-1,
\]
which is what we wanted to prove.
\end{proof}


\nin {\bf Acknowledgments:} The results in this article were obtained as an
off-shoot of a larger project, joint with Ezra Miller, whom we thank. 
We are also very grateful to Francisco Castro-Jim{\'e}nez, Jos{\'e}
Mar{\'{\i}}a 
Ucha and Mar{\'{\i}}a Isabel 
Hartillo Hermoso, who hosted us in Sevilla while we worked on this article.

%
%
%

\def\cprime{$'$}
\providecommand{\bysame}{\leavevmode\hbox to3em{\hrulefill}\thinspace}
\providecommand{\MR}{\relax\ifhmode\unskip\space\fi MR }
\providecommand{\MRhref}[2]{%
  \href{http://www.ams.org/mathscinet-getitem?mr=#1}{#2}
}
\providecommand{\href}[2]{#2}

\end{document}